\documentclass[a4paper,10pt]{article}

\usepackage[latin1]{inputenc}
\usepackage[english]{babel}
\usepackage{amsmath}
\usepackage{amsthm}
\usepackage{amssymb}
\usepackage{amsfonts}
\usepackage{graphics,graphicx}
\usepackage{subfigure}
\usepackage{url}

\DeclareSymbolFont{AMSb}{U}{msb}{m}{n}
\DeclareMathSymbol{\N}{\mathbin}{AMSb}{"4E}
\DeclareMathSymbol{\Z}{\mathbin}{AMSb}{"5A}
\DeclareMathSymbol{\R}{\mathbin}{AMSb}{"52}
\DeclareMathSymbol{\Q}{\mathbin}{AMSb}{"51}
\DeclareMathSymbol{\C}{\mathbin}{AMSb}{"43}

\newcommand{\Ma}{\mathbf{a}}
\newcommand{\Mc}{\mathbf{c}}
\newcommand{\Mu}{\mathbf{u}}

\newcommand{\Mw}{\mathbf{w}}
\newcommand{\Mx}{\mathbf{x}}

\newcommand{\MOMn}{{\mathcal M}_{n+1}}
\newcommand{\MOM}[1]{{\mathcal M}_{#1}}
\newcommand{\Pn}{{\mathcal P}_{n+1}}

\DeclareMathOperator{\SPAN}{span}
\DeclareMathOperator{\INT}{Int}

\newtheorem{lemma}{Lemma}[section]
\newtheorem{theorem}[lemma]{Theorem}

\title{On the computation of Gaussian quadrature rules for Chebyshev sets of linearly independent functions}
\author{
  Daan Huybrechs\footnote{Department of Computer Science, KU Leuven, Belgium (daan.huybrechs@cs.kuleuven.be)}
}
\date{}

\begin{document}

\maketitle
\begin{abstract}
We consider the computation of quadrature rules that are exact for a Chebyshev set of linearly independent functions on an interval $[a,b]$. A general theory of Chebyshev sets guarantees the existence of rules with a Gaussian property, in the sense that $2l$ basis functions can be integrated exactly with just $l$ points and weights. Moreover, all weights are positive and the points lie inside the interval $[a,b]$. However, the points are not the roots of an orthogonal polynomial or any other known special function as in the case of regular Gaussian quadrature. The rules are characterized by a nonlinear system of equations, and earlier numerical methods have mostly focused on finding suitable starting values for a Newton iteration to solve this system. In this paper we describe an alternative scheme that is robust and generally applicable for so-called complete Chebyshev sets. These are ordered Chebyshev sets where the first $k$ elements also form a Chebyshev set for each $k$. The points of the quadrature rule are computed one by one, increasing exactness of the rule in each step. Each step reduces to finding the unique root of a univariate and monotonic function. As such, the scheme of this paper is guaranteed to succeed. The quadrature rules are of interest for integrals with non-smooth integrands that are not well approximated by polynomials.
\end{abstract}

\section{Introduction}\label{sect:introduction}

A set $T_{n+1} := \{ u_j \}_{j=0}^n$ is called a \emph{Tchebysheff system} or a \emph{Chebyshev set}\footnote{The former name is used by Karlin and Studden in \cite{karlin1966tchebysheff} on which this paper is based, while the latter is perhaps more common in numerical analysis.} on the interval $[a,b]$ if $u_j$ are continuous real-valued functions on a closed finite interval $[a,b]$, and if furthermore the set satisfies one of the following equivalent conditions:
\begin{itemize}
 \item any real linear combination of the form $v = \sum_{j=0}^{n} c_j u_j$ has at most $n+1$ distinct zeros in the interval $[a,b]$;
 \item the determinant
 \begin{equation}\label{eq:determinant}
\left|  \begin{array}{cccc}
  u_0(t_0) & u_0(t_1) & \cdots & u_0(t_n) \\
  u_1(t_0) & u_1(t_1) & \cdots & u_1(t_n) \\
  \vdots & \vdots & & \vdots \\
  u_n(t_o) & u_n(t_1) & \cdots & u_n(t_n) \\
 \end{array} \right|
 \end{equation}
 does not vanish whenever $a \leq t_0 < t_1 < \ldots < t_n <= b$.
\end{itemize}
A classical example of a Chebyshev set are the polynomials $\{ x^j \}_{j=0}^{n}$, in which case \eqref{eq:determinant} is a \emph{Vandermonde determinant}. Eigenfunctions of Sturm-Liouville operators corresponding to the smallest eigenvalues also form a Chebyshev set. Note that the non-vanishing determinant \eqref{eq:determinant} implies, apart from linear independence of the basis functions, that the interpolation problem is uniquely solvable in any set of $n+1$ distinct points in $[a,b]$.

The study of Chebyshev sets dates back to Markov at the end of the 19th century~\cite{markov1898limitingvalues} and was subsequently extended by many others. A comprehensive theory following a geometric approach was developed by Krein \cite{krein1951} and by Karlin and Studden \cite{karlin1966tchebysheff}. In this paper we follow the geometric approach of the latter reference. We review some of its basic concepts in \S\ref{sect:review}.

The relevance of this theory to numerical quadrature was realised early on, though literature on this topic is not extensive. Karlin and Studden remark on the applicability of their results to \emph{mechanical quadrature} in \cite[\S IV.8]{karlin1966tchebysheff}, but otherwise focus on the abstract theory of moment spaces. Various results in non-polynomial Gaussian qaudrature are listed by Gautschi in a survey paper on Gauss-Christoffel quadrature formulae, in which the work of Karlin and Studden is explicitly mentioned \cite{gautschi1981christoffel}. The first numerical method for a general Chebyshev set $T_{2l} := \{ u_j \}_{j=0}^{2l-1}$ was described, to the best of our knowledge, by Ma, Rokhlin and Wandzura in \cite{ma1996generalgauss}. Their method is based on a continuation scheme for the non-linear system of equations characterizing exactness on the space spanned by the set $T_{2l}$:
\begin{equation}\label{eq:system}
 \sum_{i=1}^l \lambda_i u_j(t_i) = \int_a^b u_j(t) w(t) {\rm d}t, \qquad j=0,\ldots,2l-1.
\end{equation}
Here, $w(x) > 0$ is a positive weight function, and $t_i$ and $\lambda_i$ are the unknown points and weights of the Gaussian quadrature rule with length $l$.\footnote{We use the symbol $n$ in this paper with the same meaning as it has in the theory and notation of \cite{karlin1966tchebysheff}, i.e. the size of the Chebyshev set $T_{n}$, so that we can recall and reformulate their results without ambiguity. For this reason, we will denote the length of a quadrature rule by $l$. We aim for the correspondence $2l=n+1$, such that the generalized Gaussian quadrature rule with $l$ points is exact for the $2l$ functions in the set $T_{n+1} = T_{2l}$. Note that this corresponds to an odd value of $n$, and that the parity of $n$ matters in several results.} The equations \eqref{eq:system} define the generalized Gaussian quadrature rule and it is known that a unique solution exists for the points and weights.

Later papers by various authors include \cite{yarvin1998generalgauss,huybrechs2009generalizedgauss,bremer2010generalizedgauss,bremer2010universalquadrature}. Convergence of the quadrature rule is studied in \cite{huybrechs2009generalizedgauss} for the specific case of functions of the form
\begin{equation}\label{eq:singular_f}
f(t) = p(t) + s(t) q(t),
\end{equation}
where $p$ and $q$ are smooth functions well approximated by polynomials, and $s(t)$ is a known singularity function such as $\log t$. A simple corresponding Chebyshev set is $\{ x^j \}_{j=0}^{l-1} \cup \{ s(x) x^j \}_{j=0}^{l-1}$. One can also use Chebyshev polynomials $T_j$ instead of monomials $x^j$ without changing the span of the set and this is generally a better choice numerically.

Applications of these generalized Gaussian quadrature rules include the evaluation of integrals with endpoint singularities, of logarithmic or algebraic nature, or with various singularities that may arise in the discretization of integral equation methods with singular kernel functions \cite{bremer2010universalquadrature}. A fundamental beneficial feature of the rules is that one does not need to identify the singular and non-singular parts of the integrand -- they are incorporated into the Chebyshev set. One simply evaluates the integrand as a whole, as follows:
\begin{equation*}
 \int_a^b f(t) w(t) {\rm d}t \approx \sum_{i=1}^l \lambda_i f(t_i).
\end{equation*}
Here, even if $f$ is singular of the form \eqref{eq:singular_f}, one does not need to evaluate nor identify explicitly the smooth parts $p$ and $q$ separately. It is sufficient to be able to evaluate $f$.

In the case of Gaussian quadrature rules for regular polynomials, the quadrature points are the roots of the orthogonal polynomial of degree $l$ with respect to the weight function $w(x)$ \cite{davis1984numericalintegration,gautschi2004opq}. The orthogonality conditions of the polynomial form linear conditions on its coefficients, and hence finding an orthogonal polynomial is a linear problem. Finding its roots is a nonlinear problem, but they can be computed efficiently as the eigenvalues of the associated \emph{Jacobi matrix}, a tridiagonal matrix defined in terms of the recurrence coefficients of the orthogonal polynomial sequence \cite{gautschi2004opq}. A popular fast ${\mathcal O}(l^2)$ algorithm to do so is that of Golub and Welsch \cite{golub1969gauss}, but more recently various ${\mathcal O}(l)$ algorithms have started to appear for the classical Gaussian rules \cite{glaser2007roots,bogaert2012legendre,hale2013fastgauss}.

Orthogonal polynomials play no role in this paper. Instead, we rely on theoretical results on Chebyshev sets. Our numerical approach is to compute the quadrature rules one point at a time. We alternate between the computation of so-called \emph{principal representations} and \emph{canonical representations} of points in the moment space induced by the Chebyshev set. We review these concepts in \S\ref{sect:review} and describe the algorithm in \S\ref{sect:algorithm}. In brief, a principal representation of a $k$-dimensional moment vector corresponds to a Gaussian quadrature rule with a minimal number of points, approximately $k/2$. This representation is deformed continuously into a so-called canonical representation of the same moment vector by adding one point, precisely in one of the endpoints of $[a,b]$. A canonical representation can be thought of as a Gaussian quadrature rule with one point fixed a priori. In turn, this canonical representation is deformed continuously into a principal representation of a $k+1$-dimensional moment vector. This means that the fixed point is varied until all points agree with the next Gaussian quadrature rule of higher order. Starting from $k=1$, in each step along the way the current points and weights in the algorithm form a stable quadrature rule (with positive 
weights) that is guaranteed to exist. 
Moreover, the deformations involved are strictly monotonic and therefore straightforward to carry out numerically. In practice, the methodology means that quadrature points are added one at a time in one endpoint of the interval $[a,b]$, and each point moves monotonically in the direction of the other endpoint until it reaches its final position in the $l$-point quadrature rule. We illustrate this process with several examples.

\section{Review: a geometric theory of Chebyshev sets}\label{sect:review}

We review the concepts of the geometric theory of Chebyshev sets as formulated by Karlin and Studden, which are relevant for the formulation and understanding of our algorithm. We include a description of the link between these concepts and numerical quadrature. All the results in this section are written explicitly in  \cite{karlin1966tchebysheff}, mostly in Chapter II, and are not due to the authors of this paper. References will be given in this section only for the precise statement of theorems, yet we adopt the same notation, terminology and at times precise formulations from \cite{karlin1966tchebysheff} throughout this section.

It should be noted that we could have equally based our review on the formulation of results by Krein in \cite{krein1977momentproblem}.

\subsection{The moment space}

Consider a Chebyshev set $T_{n+1} := \{ u_j \}_{j=0}^n$ of length $n+1$ on the interval $[a,b]$. We call any linear combination of the form $\sum_{j=0}^n a_j u_j$ a \emph{generalized polynomial} or \emph{T-polynomial}. 

We define the \emph{moment space} ${\mathcal M}_{n+1}$ as
\begin{equation}\label{eq:momentspace}
 \MOMn = \{ \Mc  = (c_0,c_1,\ldots,c_n) \in \R^{n+1} \, | \, c_j = \int_a^b u_j(t) {\rm d}\sigma(t) \},
\end{equation}
where $\sigma(t)$ ranges over all nondecreasing right continuous functions of bounded variation. Thus, the moment space contains the vectors of moments of $T_{n+1}$ with respect to any possible measure. The specific measure we are interested in is ${\rm d}\sigma(t) = w(t){\rm d}t$, with $w(x) > 0$ a positive weight function.

It can be shown that the moment space ${\mathcal M}_{n+1}$ is a closed convex cone in $\R^{n+1}$. In particular, it can also be characterized as the convex conical hull of a parametric curve that is generated by the elements of the set:
\begin{equation}\label{eq:curve}
 C_{n+1} = \{ \gamma_t = (u_0(t),u_1(t),\ldots,u_n(t) ) \, | \, a \leq t \leq b \}.
\end{equation}
It is a geometric property of such cones, due to Carath{\'e}odory \cite{caratheodory1911}, that each point in such convex hull can be written as a linear combination of at most $n+2$ points on the curve:
\begin{equation}\label{eq:caratheodory}
 \forall \mathbf{\gamma} \in \MOMn: \, \, \gamma_j = \sum_{i=1}^{n+2} \lambda_i u_j(t_i), \qquad j = 0,\ldots,n,
\end{equation}
with values $\lambda_i > 0$ and $a \leq t_i \leq b$, $i=1,\ldots,n+2$.

Note the similarity between the linear combination of points in the moment space \eqref{eq:caratheodory} and the exactness conditions of a quadrature rule in \eqref{eq:system}. Indeed, the right hand side in \eqref{eq:system} for $j=0,\ldots,2l-1$, describes a moment vector in $\MOM{2l-1}$ corresponding to the known measure ${\rm d}\sigma(t) = w(t) {\rm d}t$. The left hand side is a convex (because $\lambda_i > 0$) linear combination of points on the curve $C_{2l-1}$. Vice-versa, any convex linear combination of points on the curve $C_{n+1}$ as in \eqref{eq:caratheodory} corresponds to a quadrature rule for the measure(s) associated with $\mathbf{\gamma}$, with positive weights and points inside the interval $[a,b]$, and in this case with $n+2$ points because there are $n+2$ terms in the expansion.

The most important results, for our purposes, in the theory of \cite{karlin1966tchebysheff} are the characterization of the minimal number of terms required in such linear combinations. This corresponds to quadrature rules with positive weights and with a minimal number of points. The general result of Carath{\'e}odory already guarantees the existence of quadrature rules with $n+2$ points inside the interval and with positive weights, but the generalized Gaussian quadrature rules under consideration in this paper have fewer than half as many points.

\subsection{Representations of a moment vector}

Any moment vector $\Mc$ in $\MOMn$ is finite-dimensional and can result from many different measures $\sigma$, as long as their first $n+1$ moments agree. These measures are called \emph{representations} of $\Mc$. Of particular interest are \emph{convex representations} of the form
\begin{equation}\label{eq:representation}
 \Mc = \sum_{i=1}^p \lambda_i \Mu(t_i),
\end{equation}
where we have used the notation
\[
\Mu(t) = (u_0(t), u_1(t), \ldots, u_n(t))
\]
and where $\lambda_i > 0$ and $a \leq t_i \leq b$, $i=0,\ldots,n$. We say that the representation \emph{involves} the points $t_i$, which are also called the \emph{roots} of the representation. An equivalent statement to \eqref{eq:representation} is that
\[
 c_j = \sum_{i=1}^p \lambda_i u_j(t_i), \qquad j=0,1,\ldots,n,
\]
which we recognize as having the same form as \eqref{eq:system} and \eqref{eq:caratheodory}. Each representation of the form \eqref{eq:representation} describes a quadrature rule that is exact for a certain set of moments $\Mc$. One can also think of the measure corresponding to the representation \eqref{eq:representation} as being discrete: it accumulates all its variation at the roots $t_i$, with associated weight $\lambda_i$. It happens to be the case that this discrete measure agrees with the measure ${\rm d}\sigma(t) = w(t){\rm d}t$ when restricted to $\MOMn$. The representation is called convex because $\lambda_i > 0$.

In order to count the number of points in a set, we will adopt a special rule. We define the \emph{index} of a set $\{ t_i \}_{i=1}^p$ as the number of points in the set, with the special convention that an interior point $t_i \in (a,b)$ receives a full count, but an endpoint $a$ or $b$ counts only as a half point. The index of a representation of the form \eqref{eq:representation} is the index of the set of roots it involves. Finally, since we are interested in representations with a minimal number of points, we define the index $I(\Mc)$ of a point $\Mc \in \MOMn$ as the minimal index of any of its convex representations.

We already know from \eqref{eq:caratheodory} that in general $I(\Mc) \leq n+2$. This is an upper bound, but it is not sharp: we are aiming for a value approximately half this upper bound. There are also lower bounds, and the first result we formulate expresses that points with small index must lie on the boundary of the moment cone.

\begin{theorem}[\cite{karlin1966tchebysheff}, Theorem 2.1]\label{thm:boundary}
A vector $\Mc \in \MOMn$ (with $\Mc \neq \mathbf{0}$) is a boundary point of $\MOMn$ if and only if $I(\Mc) < (n+1)/2$. Moreover, every boundary point $\Mc \in \MOMn$ admits a unique representation
\begin{equation}\label{eq:boundary_representation}
 \Mc \in \MOMn = \sum_{i=1}^p \lambda_i \Mu(t_i)
\end{equation}
where $p \leq \frac{n+2}{2}$ and $\lambda_i > 0$, $i=1,2, \ldots, p$.
\end{theorem}
This theorem means that any point in the interior of the moment cone, the general case, must have index at least $(n+1)/2$. This implies a lower bound on the number of roots, and we will see that it is sharp for any point in the interior of $\MOMn$. However, before characterizing the corresponding representations further, we first recall a result about representations with larger index, in which one root is fixed.

\begin{theorem}[\cite{karlin1966tchebysheff}, Theorem 3.1]\label{thm:canonical}
 Let $\Mc$ be an interior point of $\MOMn$. For each $t^*$, $a \leq t^* \leq b$, there exists a representation with positive weights $\lambda_i > 0$, $i=1,2,\ldots,p$,
\begin{equation}\label{eq:canonical_representation}
 \Mc \in \MOMn = \sum_{i=1}^p \lambda_i \Mu(t_i)
\end{equation}
of index $(n+1)/2$ or $(n+2)/2$ which involves the point $t^*$.
\end{theorem}

This result has interesting ramifications for numerical integration on Chebyshev sets. Indeed, one can choose any fixed point $t^* \in [a,b]$. Theorem \ref{thm:canonical} guarantees that a quadrature rule with positive weights exists with at most $(n+2)/2$ points, including $t^*$, that will integrate $n+1$ elements of the set exactly. These rules have $(n+2)/2$ points in general, but may have only $(n+1)/2$ points if we happened to choose $t^*$ as one of the roots of the Gaussian quadrature rule.

Recall that endpoints count only as a half point in this statement. Hence, if $n$ is odd, then we know that a rule with $(n+2)/2$ points has to involve exactly one of the endpoints $a$ or $b$. We will use this observation to our advantage later.

The existence of the Gaussian quadrature rule itself is related to the existence of a special so-called \emph{lower principal representation} of $\Mc$.

\subsection{Canonical and principal representations}

Let $\Mc$ be an interior point of the moment cone $\MOMn$. If a representation of $\Mc$ has index $(n+1)/2$ it is called \emph{principal}. If it has index $(n+2)/2$ it is called \emph{canonical}. The principal representations are the ones we are looking for, since they involve the smallest number of roots.

The existence of canonical representations is guaranteed by Theorem \ref{thm:canonical}. Indeed, there is at least a one-parameter family of canonical representations for any point $\Mc$, which can be described by varying the fixed point $t^*$ in the interval $[a,b]$. These are, in fact, all canonical representations of $\Mc$. For each fixed $t^*$, the canonical representation \eqref{eq:canonical_representation} is unique \cite[Corollary 3.2]{karlin1966tchebysheff}.

There are precisely two principal representations. One is called the \emph{lower principal representation} and we denote it by $\underline{\sigma}$. The other is called the \emph{upper principal representation}, denoted $\overline{\sigma}$. We will describe the points they involve in detail. This depends on the parity of $n$, but a general statement one can make is that the upper principal representation involves the right endpoint $b$. The Gaussian quadrature rule with $l$ points we are interested in corresponds, as Gautschi already points out in \cite{gautschi1981christoffel}, to the lower principal representation for odd $n=2l-1$.

The existence of the principal representations essentially follows from Theorem \ref{thm:canonical}. For even $n=2m$, it turns out one can simply choose $t^*=a$ or $t^*=b$. Since $(n+1)/2$ is a half integer in this case, and due to the way we count endpoints, each of these two representations can not involve the other endpoint. Hence, these two canonical representations are both principal. The one that involves the left endpoint is the lower principal representation, while the one with $t^*=b$ is the upper principal representation.

For odd $n = 2m+1$ it is slightly more complicated. One principal representation is found from the choice $t^*=a$, and this representation now also involves $b$ because its index $(n+1)/2$ has to be an integer. Since this representation involves the right endpoint, it is the upper principal representation. The lower principal respresentation does not involve any endpoints. Karlin and Studden describe in \cite[\S II.3]{karlin1966tchebysheff} shrinking the interval $[a,b]$ until one can again choose $t^*$ as one of the endpoints. Unfortunately, this existence proof is not constructive.

An interesting result to recall, before we characterize the lower and upper representations further, is the following interlacing property.

\begin{theorem}[\cite{karlin1966tchebysheff}, Corollary 3.1]\label{thm:principal_interlacing}
 For each $\Mc \in \INT \MOMn$ there exist precisely two principal representations. The roots of these representations strictly interlace.
\end{theorem}

It is helpful to summarize our findings on the possible distributions of roots of prinicipal representations for odd and even $n$. In the numerical method we will be using all of them, not just the lower principal representation for odd $n$.

For even $n=2m$, the principal representations have half-integer index $(n+1)/2=m+\frac12$ and we have the following set of $m+1$ roots in the lower principal representation:
\begin{equation}\label{eq:even_roots_lower}
 \underline{\sigma}: \qquad a = t_1 < t_2 < t_3 < \ldots < t_{m+1} < b \qquad \mbox{(lower)}
\end{equation}
For the upper representation, the roots include the other endpoint $b$:
\begin{equation}\label{eq:even_roots_upper}
 \overline{\sigma}: \qquad a < s_1 < s_2 < s_3 < \ldots < s_{m+1} = b \qquad \mbox{(upper)}
\end{equation}
In addition, the strict interlacing of the roots of principal representations (by Theorem \ref{thm:principal_interlacing}) implies
\begin{equation*}
 t_i < s_i < t_{i+1} < s_{i+1}, \qquad i=1,\ldots,m.
\end{equation*}

For odd $n=2m+1$, both representations have index $(n+1)/2=m+1$ but the number of roots differs since two endpoints count as one point. We have $m+1$ roots in the lower principal representation,
\begin{equation}\label{eq:odd_roots_lower}
 \underline{\sigma}: \qquad a < t_1 < t_2 < t_3 < \ldots < t_{m+1} < b \qquad \mbox{(lower)}
\end{equation}
but $m+2$ points in the upper one:
\begin{equation}\label{eq:odd_roots_upper}
 \overline{\sigma}: \qquad a = s_1 < s_2 < s_3 < \ldots < s_{m+2} = b \qquad \mbox{(upper)}
\end{equation}
The interlacing property now means that
\[
  s_i < t_i < s_{i+1}, \qquad i=1,\ldots,m+1.
\]

In our case of interest where $n +1 = 2l$, we have $m+1 = l$ and \eqref{eq:odd_roots_lower} captures the $l$ roots of a Gaussian quadrature rule exact for the $2l$-dimensional Chebyshev set $T_{n+1}=T_{2l}$.  One recognizes in \eqref{eq:odd_roots_upper} the corresponding \emph{Gauss-Lobatto rule}, the well-known variant of Gaussian quadrature that includes the endpoints. The rules \eqref{eq:even_roots_lower} and \eqref{eq:even_roots_upper} can be seen as generalizing Gauss-Radau formulas, each of them including just one endpoint \cite{davis1984numericalintegration}.

\subsection{Interlacing properties and continuity of canonical representations}

The canonical representations correspond to a fixed value $t^* \in [a,b]$ by Theorem \ref{thm:canonical}. In the algorithm of \S\ref{sect:algorithm} we will be computing canonical representations precisely by varying $t^*$, hence we include more theory.

We have seen the interlacing properties of principal representations in Theorem \ref{thm:principal_interlacing}. In fact, a much stronger statement holds that covers all canonical representations of a moment vector $\Mc$: the roots of any two canonical representations must interlace. This leads to strong restrictions on where those roots can be.

\begin{theorem}[\cite{karlin1966tchebysheff}, Theorem 3.2]\label{thm:canonical_interlacing}
 Let $\Mc \in \INT \MOMn$.
 
Consider two different representations $\sigma'$ and $\sigma''$ of $\Mc$ with index $\leq (n+2)/2$ and with roots $\{t_i'\}_1^p$ and $\{t_i''\}_1^q$, and weights $\{\lambda_i'\}_1^p$ and $\{\lambda_i''\}_1^q$. Then the roots $\{t_i'\}_1^p$ and $\{t_i''\}_1^q$ strictly interlace in the open interval $(a,b)$, but they may possibly share one or both of the endpoints $a$ or $b$.
 
Moreover, if $t_1' = t_1'' = a$ then $\lambda_1' \neq \lambda_1''$ and $\lambda_1' > \lambda_1''$ if and only if $t_2' > t_2''$. Similarly, if $t_p' = t_q'' = b$ then $\lambda_p' \neq \lambda_q''$ and $\lambda_p' > \lambda_q''$ if and only if $t_{p-1}' < t_{q-1}''$.
\end{theorem}
The second half of the theorem expresses that roots at an endpoint can coincide, but the corresponding weight is greater for the representation whose neighbouring root is farthest away.

\subsubsection{Odd $n = 2m+1$}

Consider an odd value of $n=2m+1$ and recall the distribution of the roots of the lower and upper principal representations in \eqref{eq:odd_roots_lower} and \eqref{eq:odd_roots_upper}. We define two types of intervals between two consecutive points of the different representations:
\begin{align}\label{eq:canonical_intervals_odd}
 K_i &= (s_i, t_i), \qquad i=1,2,\ldots,m+1, \\
 J_i &= (t_i, s_{i+1}), \qquad i=1,2,\ldots,m+1. \nonumber
\end{align}
For any fixed point $t^*=\xi$, the associated canonical representation has index $(n+2)/2 = m+3/2$ and it must involve $m+2$ roots including one endpoint. We will denote the points of the canonical representation by $\{ t_i^*(\xi) \}_{i=1}^{m+2}$ as a function of $\xi$. Now let the point $\xi$ vary in $K_1$, i.e. from the left endpoint $s_1=a$ to the first root $t_1$ of $\underline{\sigma}$. Due to the interlacing property, it must be the case that the roots of the canonical representation satisfy $t_i^*(\xi) \in K_i$ for $i=1,2,\ldots,m+1$ and $t_{m+2}^*(\xi) = b$. That is, we must have that
\begin{align}
a = s_1 < \xi = t_1^* < t_1 &< s_2 < t_2^* < t_2 < s_3 < t_3^* < \ldots \label{eq:canonical_odd_K} \\
& \ldots < s_{m+1} < t_{m+1}^* < t_{m+1} < s_{m+2} = t_{m+2}^* = b. \nonumber
\end{align}

On the other hand, when $\xi$ varies further in $J_1 = (t_1,s_2)$, the canonical representation with $m+2$ roots must include the left endpoint $a$. We now label the points as $\{t_i^*(\xi)\}_0^{m+1}$, with $t_0^*=a$. The interlacing properties prescribe that $t_i^*(\xi) \in J_i$, $i=1,\ldots,m+1$:
\begin{align}
a = s_1 = t_0^* < t_1 < \xi = t_1^* &< s_2 < t_2 < t_2^* < s_3 < t_3 < \ldots \label{eq:canonical_odd_J} \\
& \ldots < s_{m+1} < t_{m+1} < t_{m+1}^* <  s_{m+2} = b. \nonumber
\end{align}
With our choice of labelling, it is as if for $\xi > t_1$ the point $t^*_{m+2}$ leaves the interval to the right, whereas $t_0^*$ enters from the left.

In both cases, the points $t_i(\xi)$ monotonically traverse their allowed intervals from left to right. Indeed they have to: any lack of monotonicity would lead to a violation of the interlacing property of two canonical represenations nearby. Moreover, all points are continuous functions of $\xi$.

It remains to describe what happens in the corner cases $\xi=a$, $\xi=t_1$ and $\xi=s_1$. As it turns out, the limits to all these cases are free of singularities. At $\xi=0$, the weight $\lambda_{m+2}^*(\xi)$ at $t_{m+2}^*=b$ is maximal and it decreases monotonically to $0$ at $\xi=t_1$, after which the point $b$ is no longer included. The monotonic decrease of this weight is due to the second part of Theorem \ref{thm:canonical_interlacing} and the fact that $t_{m+1}^*(\xi)$ moves monotonically closer to $t_{m+1}$ (and thus to $b$). Similarly, for the value of $\xi=t_1$ and onwards, the weight $\lambda_0^*(\xi)$ at $t_0^*(\xi)=a$ increases monotonically from $0$ to its maximal value at $\xi=s_1$. Thus, at the value $\xi = t_1$, both $a$ and $b$ appear to be included simultaneously, but they both have weight $0$. The canonical representation is in fact the lower principal representation \eqref{eq:odd_roots_lower}. For smaller $\xi$ the right endpoint is involved, for larger $\xi$ the left endpoint is, but 
the transition is at least continuous.

Due to the uniqueness of canonical representations, it is sufficient to traverse the intervals $J_1$ and $K_1$. As $\xi$ progresses further, no new canonical representations are encountered: all intervals $J_i$ and $K_i$ were traversed monotonically from left to right already by $t_i^*(\xi)$ along with $\xi$. This process simply repeats itself for larger $\xi$.

\subsubsection{Even $n = 2m$}

We also include a description of what happens for even $n=2m$. Recall the lower and upper principal representations \eqref{eq:even_roots_lower} and \eqref{eq:even_roots_upper}. In this case, we define the intervals
\begin{align}\label{eq:canonical_intervals_even}
 J_i &= (t_i, s_i), \qquad i=1,2,\ldots,m+1, \\
 K_i &= (s_i, t_{i+1}), \qquad i=1,2,\ldots,m. \nonumber
\end{align}

As $\xi$ traverses the interval $J_1$, i.e. from $t_1=a$ to $s_1$, we find due to the interlacing properties that there are $m+1$ roots in the canonical representation satisfying $t_i^*(\xi) \in J_i$, $1 \leq i \leq m+1$. More specifically, we have
\begin{align}
a = t_1 < \xi = t_1^* < s_1 &< t_2 < t_2^* < s_2 < t_3 < t_3^* < \ldots \label{eq:canonical_even_J}\\
& \ldots < t_{m+1} < t_{m+1}^* < s_{m+1} = b. \nonumber
\end{align}
Continuing with $\xi$ traversing $K_1$, from $s_1$ to $t_2$, we find that there must now be $m+2$ roots, including both endpoints. We have $t_0^*(\xi) = a$, $t_1^*(\xi)=\xi$, $t_{m+1}^*(\xi)=b$ and $t_i^*(\xi) \in K_i$, $i=1,\ldots,m$. Summarizing:
\begin{align}
a = t_0^* = t_1 <  s_1 < \xi= t_1^* &< t_2 < s_2 < t_2^* < t_3 < s_3 < \ldots \label{eq:canonical_even_K} \\
& \ldots < s_m < t_m^* < t_{m+1} < s_{m+1} = t_{m+1}^* = b.\nonumber
\end{align}
At the transition $\xi=s_1$, both endpoints become involved: $a$ enters with weight $0$, while $b$ is approached by $t_{m+1}^*(\xi)$ and has maximal weight.  As $\xi$ traverses $K_1$, the weight of $a$ increases to its maximal value, while the weight of $b$ decreases down to $0$. As in the case of odd $n$, there is no discontuinity in the process, and the process simply repeats itself as $\xi$ traverses $J_2$ and $K_2$ with no new canonical representations encountered.

\subsection{The range of a moment vector}
\label{ss:range}

The names of the lower and upper principal representations originate in a property called the \emph{range} of a moment vector. Assume that $T_{n+2}$ is also a Chebyshev set, with one additional element $u_{n+2}$ compared to $T_{n+1}$. The range $R(\Mc)$ is the set of values
\begin{equation}\label{eq:range}
 \int_a^b u_{n+1}(t) {\rm d}\sigma(t),
\end{equation}
where $\sigma$ varies over all measures representing $\Mc$, i.e. over all measures that have the same moment vector $\Mc$ in $\MOMn$.

The range is a closed interval that can be described by its minimal and maximal values:
\[
 R(\Mc) = \{ \gamma : \underline{\gamma} < \gamma < \overline{\gamma} \}.
\]
The moment $\underline{\Mc} = (c_0, c_1, \ldots, c_n, \underline{\gamma})$, an extension of $\Mc$ to dimension $n+2$, lies on the boundary of the $n+2$-dimensional moment cone $\MOM{n+2}$. So does $\overline{\Mc} = (c_0, c_1, \ldots, c_n, \overline{\gamma})$. Therefore, by Theorem \ref{thm:boundary}, these moment vectors have an index less than or equal to $(n+1)/2$. It can be shown that their index is, in fact, precisely $(n+1)/2$, and that the corresponding representations are precisely the lower and upper principal representations of $\Mc$ in $\MOMn$.

This property explains the names of the principal representations. The lower principal representation of $\Mc$ corresponds to the lower point $\underline{\gamma}$ of the range $R(\Mc)$, while the upper principal representation corresponds to the upper point $\overline{\gamma}$.

For a given moment vector $\Mc \in \MOMn$, we can consider all intermediate moment vectors in $\MOM{n+2}$ of the form
\[
\Mc_\gamma = (c_0, c_1, \ldots, c_n, \gamma)
\]
with $\underline{\gamma} \leq \gamma \leq \overline{\gamma}$. These vectors are all points on the line connecting the two boundary points $\underline{\Mc}$ and $\overline{\Mc}$ of $\MOM{n+2}$. Each such vector has two distinct principal representations, of index $(n+2)/2$. Hence, these are also canonical representations of the original moment vector $\Mc$. The two principal representations of $\Mc_\gamma$ form another parameterization of the canonical representations of $\Mc$.

The last result we shall need for our algorithm is the following description of the map onto the range of $\Mc$ for varying $\gamma$. The \emph{upper and lower canonical intervals} in the formulation of this theorem are the intervals $J_i$ and $K_i$ respectively, defined by \eqref{eq:canonical_intervals_odd} or \eqref{eq:canonical_intervals_even}. The $J_i$'s are called upper because they include the right endpoint $b$.

\begin{theorem}[\cite{krein1951,karlin1953geometry}]\label{thm:range_map}
 Let $\Mc \in \INT \MOMn$ and let $I$ be the closure of one of the upper or lower canonical intervals. Then the canonical representation $\sigma_\xi$ (principal of $\xi=a$ or $b$) induce a continuous $1:1$ mapping $\gamma(\xi)$, of the interval $I$ onto the interval $R(\Mc)$, defined by
 \begin{equation}\label{eq:range_map}
  \gamma(\xi) = \int_a^b u_{n+1}(t) {\rm d}\sigma_\xi(t).
 \end{equation}
 The function $\gamma(\xi)$ is increasing on $I = \overline{J}_i$ and decreasing on $I = \overline{K}_i$.
\end{theorem}

Of particular interest to us is the monotonicity of this map. By varying $\xi$ and computing the canonical representations of $\Mc$, we can generate principal representations of a higher-dimensional moment vector, where we can force the last (and new) entry to be anything in its allowed range. This is the main tool upon which our algorithm is based: we match the moments in the right hand side of \eqref{eq:system} one by one, each time lifting a canonical representation of $\Mc \in \MOM{k}$ to a principal representation of $\Mc_\gamma \in \MOM{k+1}$. The fact that the map is monotonic will make the numerical continuation procedure particularly straightforward.

\subsection{Nonnegative polynomials}
\label{ss:nonnegative}

The theory of nonnegative polynomials plays a major role in \cite{karlin1966tchebysheff}, but not in our algorithm -- at least not currently. For the completeness of this review, and for the purpose of possible future developments, we do include a brief description.

A generalized polynomial $p = \sum_{j=0}^n a_j u_j$ is nonnegative if $p(x) \geq 0$ on $[a,b]$. There is a duality between the space of nonnegative polynomials $\Pn$ and the moment cone $\MOMn$. The dual $\zeta^+$ of a convex cone $\zeta \subset \R^{n+1}$ is defined by
\begin{equation}\label{eq:dualcone}
 \zeta^+ = \{ \Ma \in \R^{n+1} \, | \, (\Ma,\Mc) \geq 0 \mbox{~for all~} \Mc \in \zeta \}
\end{equation}
where $(\Ma,\Mc) = \sum_{j=0}^n a_j c_j$. With this notation, it is shown in \cite{karlin1966tchebysheff} that $\MOMn^+ = \Pn$ and $\Pn^+ = \MOMn$. Indeed, for any $\Ma \in \MOMn^+$ we have by definition that $\sum_{j=0}^n a_j c_j \geq 0$ for all $c \in \MOMn$, and this implies in particular that $u(t) = \sum_{j=0}^n a_j u_j(t) \geq 0$ for $t \in [a,b]$, hence $u \in \Pn$.

In the case of regular polynomials, one may think of nonnegative polynomials as having only double roots in the interval $(a,b)$, i.e. points where the polynomial vanishes along with its derivative such that there is no change of sign. (For completeness, note that nonnegative polynomials can have simple roots in the endpoints, and this observation relates to the special counting rule for endpoints.) It does not always make sense to speak of roots with higher multiplicity in general Chebyshev sets, since the elements $u_j$ of the set are not necessarily differentiable. Nonnegative polynomials with roots can be described with a limiting procedure.

A nonnegative polynomial can always be found for any Chebyshev set $T_n$ with zeros at approximately $n/2$ prescribed points -- approximately because the precise result depends on the parity of $n$ and whether or not endpoints are included. This result is shown in Chapter I, Theorems 5.1 and 5.2 of \cite{karlin1966tchebysheff}, and many of the subsequent developments are based on this property. One can think of the existence of such nonnegative polynomials as a substitute for the fact that zeros can be factored out in regular polynomials.

\subsection{Differences to the theory of orthogonal polynomials}

The theory of Gaussian quadrature has predominantly focused on orthogonal polynomials. However, the link to the concept of orthogonality exists only for regular polynomials, since it is crucially based on the fact that polynomials can be factored. Generalized polynomials over Chebyshev sets can not in general be factored. Indeed, generalized polynomials are not necessarily closed under multiplication.

The link between orthogonality and the factoring of polynomials is clear in the following reasoning. Imagine a set of $l$ quadrature points $\{x_i\}_{i=1}^l$, and consider a non-trivial regular polynomial $q$ of degree $l \leq m < 2l$ that vanishes at all these points. The quadrature approximation to the integral of this polynomial is zero, since $q$ vanishes by construction at each quadrature point. Hence, in order to retain exactness, the integral of the polynomial itself has to vanish as well. This restriction we can only satisfy by suitably choosing the quadrature points. Since polynomials can be factored, we can write $q = p r$ as a multiple of the monic polynomial $p$ of degree $l$ that vanishes at the same $l$ points and another polynomial $r$ of degree less than $l$. The condition that the integral of $q$ must vanish becomes
\begin{equation}\label{eq:orthogonality}
\int_a^b q(x) w(x) {\rm d}x = \int_a^b p(x) r(x) w(x) {\rm d}x = 0.
\end{equation}
This condition will hold for all $r$ up to degree $l-1$ if and only if $p$ is orthogonal to all lower degree polynomials with respect to $w(x)$.

Generalizing this description, one may want to characterize the subspace of all functions in the span of $T_{2l}$ that vanish at a set of $l$ points. Subsequently, in order to obtain a Gaussian quadrature rule, one wants to ensure that the integrals of all these functions vanish (for a formal statement of this condition, see \cite[Theorem 3.1]{cools1997cubature}). This reasoning was explored for a particular type of Chebyshev sets in \cite{huybrechs2009generalizedgauss} and the quadrature points could be identified as the roots of a certain generalized polynomial, but that polynomial is not characterized by orthogonality nor by any other set of linear conditions.

Note that the square of the regular orthogonal polynomial above, $p^2$, is a nonnegative polynomial on $[a,b]$. The role of orthogonal polynomials in Gaussian quadrature is, if anything, replaced by nonnegative polynomials for generalized Gaussian quadrature. This nonnegative polynomial can be characterized by an extremal property, associated with the fact that the lower principal representation achieves the lower bound in the range of a moment vector as described in \S\ref{ss:range}. As an alternative to our algorithm below, this extremal property could be the basis for a computational procedure. A recent method for the computation of regular Gaussian quadrature using methods of optimization was described in \cite{ryu2015gausslp}. As is alluded to in the paper \cite{ryu2015gausslp} itself, this method may also apply to generalized Gaussian quadrature.

\section{An algorithm for generalized Gaussian quadrature}
\label{sect:algorithm}

The general strategy of our algorithm is recursion. Starting from a Gaussian quadrature rule of length $k$, a sequence of four steps results in a Gaussian quadrature rule of length $k+1$. The main reason for repeating some results of \cite{karlin1966tchebysheff} in the previous section has been to demonstrate that this sequence of steps is guaranteed to converge and that the continuations are, in fact, simple to carry out numerically due to the monotonicity properties.

Once the steps are identified in \S\ref{ss:algorithm1}, we formulate the same algorithm in a different way that is amenable to a straightforward implementation in \S\ref{ss:algorithm2}.

\subsection{Algorithm 1}
\label{ss:algorithm1}

Define the finite moment vector $\Mc^{n+1} \in \MOMn$ for our measure ${\rm d}\sigma(t) = w(t) {\rm d}t$, i.e.:
\[
 c^{n+1}_j = \int_a^b u_j(t) w(t) {\rm d}t, \qquad j=0,\ldots,n.
\]
The lower principal representation of $\Mc^{2k}$ has index $(n+1)/2 = k$, since $n=2k-1$. Hence, it involves $k$ roots and weights, which we denote by $\{t_{2k,i}\}_1^k$ and $\{\lambda_{2k,i} \}_1^k$ respectively. We assume that we know these or, in other words, we assume that we have already computed the Gaussian quadrature rule exact on the span of $T_{2k}$.

The four steps to compute the $k+1$-rule are:

\begin{enumerate}
\item We add the right endpoint $b$ to the existing set of points, with associated weight $0$. The result is a canonical representation of $\Mc^{2k}$ with index $k+1/2$. To be precise, we define the new points $\{ \tilde{t}_{2k,i} \}_1^{k+1}$ and weights $\{ \tilde{\lambda}_{2k,i} \}_1^{k+1}$ by:
\begin{align*}
 &\tilde{t}_{2k,i} = t_{2k,i}, \qquad i=1,2,\ldots,k,
 &\tilde{t}_{2k,k+1} = b,
\end{align*}
and
\begin{align*}
 &\tilde{\lambda}_{2k,i} = \lambda_{2k,i}, \qquad i=1,2,\ldots,k.
 &\tilde{\lambda}_{2k,k+1} = 0.
\end{align*}
This step involves no computation.

\item We now consider the points $\tilde{t}_{2k,i}(\xi)$ as functions of $\xi$, and identify the case of step $1$ with the value $\xi = t_{2k,1}$. We let $\xi$ decrease from $t_{2k,1}$ down to $a$, such that all points $\tilde{t}_{2k,i}(\xi)$ are decreasing functions of $\xi$, with the exception of $\tilde{t}_{2k,k+1}(\xi)=b$ which remains constant for this range of $\xi$. We monitor the value of the next moment of our discrete measure as a function of $\xi$:
\begin{equation} \label{eq:step2_discretemoment}
 \mu_{2k}(\xi) = \sum_{i=1}^{k+1} \tilde{\lambda}_{2k,i}(\xi) u_{2k}(\tilde{t}_{2k,i}(\xi)).
\end{equation}
This function is monotonically decreasing from its maximal value (among all possible measures with moment vector $\Mc^{2k} \in \MOM{2k}$) at $\xi = t_{2k,1}$ down to its minimal value at $\xi=a$. By continuation we locate numerically the unique value $\xi^*$ for which it equals the next continuous moment:
\begin{equation}\label{eq:optimal_xi_step2}
\mu_{2k}(\xi^*) = \int_a^b u_{2k}(t) w(t) {\rm d}t.
\end{equation}
The roots $\{ \tilde{t}_{2k,i}(\xi^*) \}_1^{k+1}$ and weights $\{ \tilde{\lambda}_{2k,i}(\xi^*)\}_1^{k+1}$ correspond to the upper principal representation of $\Mc^{2k+1}$ with index $k+1/2$ as described by \eqref{eq:even_roots_upper} for even $n=2k$.

\item We define a new set of points $\{ t_{2k+1,i} \}_1^{k+1}$ and weights $\{ \lambda_{2k+1,i} \}_1^{k+1}$ by
\begin{align*}
 &t_{2k+1,i} = \tilde{t}_{2k,i}(\xi^*), \qquad i=1,2,\ldots,k+1, \\
 &\lambda_{2k+1,i} = \tilde{\lambda}_{2k,i}(\xi^*), \qquad i=1,2,\ldots,k+1.
\end{align*}
This is still the upper principal representation of $\Mc^{2k+1}$, and this step requires no computation.

\item We now consider the points $t_{2k+1,i}(\xi)$ as functions of $\xi$ and identify the case of step $3$ with the value $\xi = t_{2k+1,1}$. For this value of $\xi$, the point set has index $k+1/2$. Similar to step 2, we let $\xi$ decrease from $t_{2k+1,1}$ down to $a$, such that all points $t_{2k+1,i}(\xi)$ are decreasing functions of $\xi$. The point $t_{2k+1,k+1}(\xi)$ decreases too in this step, away from $b$. Hence, for any $\xi < t_{2k+1,1}$ the points have index $k+1$. We monitor the value of the next moment of our discrete measure as a function of $\xi$:
\begin{equation} \label{eq:step4_discretemoment}
 \mu_{2k+1}(\xi) = \sum_{i=1}^{k+1} \lambda_{2k+1,i}(\xi) u_{2k+1}(t_{2k+1,i}(\xi)).
\end{equation}
By continuation we locate numerically the unique value $\hat{\xi}$ for which it equals the next continuous moment:
\begin{equation}\label{eq:optimal_xi_step4}
\mu_{2k+1}(\hat{\xi}) = \int_a^b u_{2k+1}(t) w(t) {\rm d}t.
\end{equation}
The roots $\{ t_{2k+1,i}(\hat{\xi}) \}_1^{k+1}$ and weights $\{ \lambda_{2k+1,i}(\hat{\xi})\}_1^{k+1}$ correspond to the lower principal representation of $\Mc^{2k+2}$ with index $k+1/2$ as described by \eqref{eq:odd_roots_lower} for odd $n=2k+1$.

\end{enumerate}

In this version of the algorithm, new points are added in the right endpoint $b$ and all points move continuously from right to left. An intermediate result is the computation of the upper principal representation of $\Mc^{2k+1}$.

Variants of this algorithm can be formulated to compute any lower or upper principal representation for any moment vector $\Mc^k$, based on the theory in the previous section. For example, we could have started in step 1 by adding the left endpoint $a$ with weight $0$, and let $\xi$ increase from $t_1$ to $s_2$. This would lead to the points moving from left to right in the continuations. We have made the choice of adding $b$ mostly because we had in mind considering integrals with singularities in the left endpoint in our examples. Adding the left endpoint to the quadrature rule would in this case lead to numerical issues.

\subsection{Discussion of the algorithm}

We make some additional comments about Algorithm 1, in order to better highlight its connection to the theory in \S\ref{sect:review}.

In step 2 of the algorithm, we start with a canonical representation of $\Mc^{2k}$. The fixed root $\xi$ just happens to be equal to the first root of the lower principal representation of $\Mc^{2k}$. In this step, $n = 2k-1$ is odd, and the interlacing properties are given by \eqref{eq:canonical_odd_K}. The fixed root traverses the interval $K_1 = (s_1,t_1)$, as defined in \eqref{eq:canonical_intervals_odd}. We start from $t_1$, the first root of the lower principal representation of $\Mc^{2k}$, and move in the direction of $s_1=a$, the first root of the upper principal representation of $\Mc^{2k}$. By Theorem~\ref{thm:range_map}, the discrete moment \eqref{eq:step2_discretemoment} is strictly decreasing on $K_1$ from its maximal value to its minimal value. (Of course, it is increasing if we traverse the interval in the opposite direction.) Hence, we are guaranteed to encounter a solution to the problem in \eqref{eq:optimal_xi_step2}, since the right hand side of \eqref{eq:optimal_xi_step2} is necessarily 
somewhere in the possible range. Note that $u_{2k}$ in \eqref{eq:optimal_xi_step2} corresponds exactly to $u_{n+1}$ in \eqref{eq:range_map}.

In step 4 of the algorithm, we start with a canonical representation of $\Mc^{2k+1}$, that is also its upper principal representation. In this step $n=2k$ is even, and the interlacing properties are given by \eqref{eq:canonical_even_J}. The fixed root now traverses the interval $J_1 = (t_1,s_1)$ as defined by \eqref{eq:canonical_intervals_even}. We start at $s_1$, the first root of the upper principal representation of $\Mc^{2k+1}$, and move again in the direction of $a = t_1$. This time, the discrete moment \eqref{eq:step4_discretemoment} is strictly increasing on $J_1$ by Theorem~\ref{thm:range_map}, from its minimal value to its maximal value. By this monotonicity, and since the right hand side of \eqref{eq:optimal_xi_step4} is necessarily somewhere in the allowed moment range, also problem \eqref{eq:optimal_xi_step4} is guaranteed to have a solution that we will encounter.

\subsection{A reformulation of the algorithm}
\label{ss:reformulation}

We intend to reformulate Algorithm 1 in a way that is more amenable to implementation and that better illustrates the computational cost. The reformulation is based on two observations. First, each quadrature rule of length less than $l$ can be extended to a quadrature rule of length $l$ by adding points with an associated weight that is zero. We can choose to put all those extra points at the right endpoint $b$. Second, the parameter $\xi$ in Algorithm 1 is actually continuous across the different steps of the algorithm. The starting value of $\xi$ in step 4 is equal to the minimal value $\xi^*$ of step 2, and the starting value of $\xi$ in step 2 is equal to the minimal value $\hat{\xi}$ in step 4 for the previous value of $k$. In fact, $\xi$ is always equal to the first root.

This means that we can perform continuation on a single global parameter, again called $\xi$, on two vectors of $l$ points $\Mx^l(\xi)$ and $l$ weights $\Mw^l(\xi)$. We have that $x_1(\xi) = \xi$. All points and weights are continuous functions of $\xi$, but they are only piecewise smooth. The breakpoints are those values of $\xi$ where an additional quadrature point comes into play.

In order to be precise, we establish some notation:
\begin{itemize}
\item We denote the vector that consists of the first $k$ components of $\Mx$ by $\Mx^k$, and we similarly use the notation $\Mw^k$. The vector components are $x_i$ and $w_i$, for $i=1,\ldots,k$. (Note that we had already defined the finite moment vector $\Mc^k$ that contains the first $k$ moments, but in that context we start counting from zero in accordance with the notation of \cite{karlin1966tchebysheff}: the elements of $\Mc^k$ are $c_j$, $j=0,\ldots,k-1$.)

\item Let $\underline{\xi}_k$ be the first root of the lower principal representation of the moment vector $\Mc^{2k} \in \R^{2k}$, and let $\overline{\xi}_k$ be the first root of the upper principal representation of the moment vector $\Mc^{2k+1} \in \R^{2k+1}$.

\item Let $\Mw = W^k(\Mx)$ be the linear map from $k$ points $\Mx \in \R^k$ to $k$ weights $\Mw \in \R^k$ such that
\[
 \sum_{i=1}^k w_i u_j(x_i) = c_j, \qquad j=0,\ldots,k-1.
\]
In other words, these are the weights of the interpolatory quadrature rule in these $k$ points.

\item Furthermore, let $\overline{X}^{k+1}(\xi)$ be the map from $\xi$ to the unique ordered set of $k+1$ roots involved in the canonical representation of the moment vector $\Mc^{2k}$, with each individual root given by $\overline{X}^{k+1}_i(\xi)$, $i=1,\ldots,k+1$. Note that $\overline{X}^{k+1}_{1}(\xi) = \xi$ and $\overline{X}^{k+1}_{k+1}(\xi) = b$. 

\item Finally, let $\underline{X}^{k+1}(\xi)$ be the map from $\xi$ to the unique ordered set of $k+1$ roots involved in the canonical representation of the moment vector $\Mc^{2k+1}$, with each individual root given by $\underline{X}^{k+1}_i(\xi)$, $i=1,\ldots,k+1$. Note that $\underline{X}^{k+1}_{1}(\xi) = \xi$ and $\underline{X}^{k+1}_{k+1}(\xi) \leq b$.
\end{itemize}

The breakpoints are aligned as follows:
\[
a <  \underline{\xi}_l < \overline{\xi}_{l-1} <  \underline{\xi}_{l-1} < \ldots < \overline{\xi}_1 < \underline{\xi_1} < \overline{\xi}_0 = b.
\]
On each interval of the form $[\overline{\xi}_k < \underline{\xi_k}]$ we are using $k+1$ points and weights and performing continuation on canonical representations of $\Mc^{2k}$ in order to find the upper principal representation of $\Mc^{2k+1}$. Subsequently, on intervals of the form $[ \underline{\xi}_{k+1}, \overline{\xi}_k]$ we are still using $k+1$ points and weights but performing continuation on canonical representations of $\Mc^{2k+1}$ in order to find the lower principal representation of $\Mc^{2k+2}$.

With our notation, this means that the continuous functions $w_i(\xi)$ and $x_i(\xi)$ satisfy the following properties:
\[
w_i(\xi) = \left\{ 
 \begin{array}{ll}
 0, & \mbox{if~} \xi > \underline{\xi}_{i-1}, \\
 W^k_i(\Mx^k), & \mbox{otherwise}.
 \end{array}\right.
\]
and
\[
x_i(\xi) = \left\{ 
 \begin{array}{ll}
 b, & \mbox{if~} \xi > \underline{\xi}_{i-1}, \\
 \overline{X}^k_i(\Mx^k), & \mbox{if~} \xi \in [\overline{\xi}_k < \underline{\xi_k}],\\
 \underline{X}^k_i(\Mx^k), & \mbox{if~} \xi \in [\underline{\xi}_{k+1} < \overline{\xi_k}].
 \end{array}\right.
\]

\subsection{Algorithm 2}
\label{ss:algorithm2}

We have to determine the breakpoints $\underline{\xi}_k$ and $\overline{\xi}_k$, as well as the associated quadrature points.

\begin{enumerate}
 \item Set $\xi = b$, $\overline{\xi}_0 = b$ and initialize $x_i(b) = b$, $i=1,\ldots,l$, $w_1(b) = c_0$ and $w_i(b) = 0$, $i = 2,\ldots,l$.
 \item Starting from $\overline{\xi}_0$ and decreasing $\xi$, solve for $\xi$ the equation
 \begin{equation}\label{eq:firstproblem}
  w_1(\xi) u_1(x_1(\xi)) =   w_1(\xi) u_1(\xi) = c_1
 \end{equation}
 while maintaining that $\Mw^1(\xi) = W^1(\Mx^1(\xi))$. Denote the solution by $\underline{\xi}_1$.
 \item Repeat for $k = 1,\ldots,l-1$:
  \begin{enumerate}
   \item Starting from $\underline{\xi}_k$ and decreasing $\xi$, solve for $\xi$ the equation
 \begin{equation}\label{eq:secondproblem}
 F_k(\xi) = \sum_{i=1}^{k+1} w_i(\xi) u_{2k}(x_i(\xi)) - c_{2k} = 0,
 \end{equation}
 while maintaining that $\Mw^{k+1}(\xi) = W^{k+1}(\Mx^{k+1}(\xi))$ and $\Mx^{k+1}(\xi) = \overline{X}^{k+1}(\xi)$. Denote the solution by $\overline{\xi}_k$.
   \item Starting from $\overline{\xi}_k$ and decreasing $\xi$, solve for $\xi$ the equation
 \begin{equation}\label{eq:thirdproblem}
 G_k(\xi) = \sum_{i=1}^{k+1} w_i(\xi) u_{2k+1}(x_i(\xi)) - c_{2k+1} = 0,
 \end{equation}
 while maintaining that $\Mw^{k+1} = W^{k+1}(\Mx^{k+1})$ and $\Mx^{k+1}(\xi) = \underline{X}^{k+1}(\xi)$. Denote the solution by $\underline{\xi}_{k+1}$.
 \end{enumerate}
\end{enumerate}

Strictly speaking, the above is not an algorithm. Rather, it is an algorithmic description of the piecewise smooth functions $w_i(\xi)$ and $x_i(\xi)$. It is not specified how these functions should be computed or stored. The formulation `for decreasing $\xi$, solve' does suggest a continuation procedure. We discuss one possible implementation next.

\subsection{Implementation}
\label{ss:implementation}

We elaborate on the solvers for the problems in the second version of the algorithm. In the first problem, given by \eqref{eq:firstproblem}, the constraint $\Mw^1 = W^1(\Mx^1)$ is simply
\[
w_1(\xi) u_0(\xi) = c_0,
\]
hence $w_1(\xi) = c_0 / u_0(\xi)$. Equation \eqref{eq:firstproblem} becomes
\[
\frac{u_1(\xi) }{u_0(\xi)} = \frac{c_1}{c_0}.
\]
The solution is the $1$-point generalized Gaussian quadrature rule that is exact for $u_0$ and $u_1$. This rule can in most cases be found explicitly.

In problems \eqref{eq:secondproblem} and \eqref{eq:thirdproblem}, the constraint on the weights is simply a linear system of equations. On the other hand, the condition on the points is a nonlinear problem. For \eqref{eq:secondproblem}, the points and weights can be solved for simultaneously from the exactness conditions:
\begin{equation}\label{eq:exact1}
 \sum_{i=1}^{k+1} w_i u_j(x_i) = c_j, \qquad j = 0, \ldots, 2k-1.
\end{equation}
There are 2k+2 points and weights in total. However, the two points $\xi_1(\xi)= \xi$ and $\xi_{k+1}(\xi) = b$ are fixed. Hence, there are $2k$ unknowns and this matches the number of equations. An assumption here is that, since we are stepping in $\xi$, the values of the previous step are good starting values for Newton's method on this nonlinear system.

The exactness conditions for the points and weights in the other problem \eqref{eq:thirdproblem} are very similar, but with one more equation:
\begin{equation}\label{eq:exact2}
 \sum_{i=1}^{k+1} w_i u_j(x_i) = c_j, \qquad j = 0, \ldots, 2k.
\end{equation}
Since $x_{k+1}(\xi)$ is no longer fixed, there is also one more degree of freedom.

Ultimately, both problems \eqref{eq:secondproblem} and \eqref{eq:thirdproblem} are about finding the root of the univariate functions $F_k$ and $G_k$, which are either motonically increasing or decreasing. We perform Newton's method again. Differentiating the left hand side of \eqref{eq:secondproblem} we obtain:
\begin{equation}\label{eq:secondproblem_derivative}
 F_k'(\xi) = \sum_{i=1}^{k+1} w_i'(\xi) u_{2k}(x_i(\xi)) + w_1(\xi)u_{2k}'(\xi) + \sum_{i=2}^k w_i(\xi) u_{2k}'(x_i(\xi)) x_i'(\xi),
\end{equation}
where the range of the second summation is smaller because we know that $x_1(\xi) = \xi$ and $x_{k+1}(\xi)=b$, hence $x_1'=1$ and $x_{k+1}'=0$. The derivatives of the other points and weights, as functions of $\xi$, can be found from the Jacobian of the nonlinear system of exactness conditions \eqref{eq:exact1}. This leads to having to solve the linear system
\begin{equation}\label{eq:exact1_derivative}
 \sum_{i=1}^{k+1} w_i' u_j(x_i) + w_1 u_j'(\xi) + \sum_{i=2}^k w_i u_j'(x_i) x_i' = 0, \qquad j = 0, \ldots, 2k-1.
\end{equation}

Similarly, differentiating the left hand side of \eqref{eq:thirdproblem} we find
\begin{equation}\label{eq:thirdproblem_derivative}
 G_k'(\xi) = \sum_{i=1}^{k+1} w_i'(\xi) u_{2k+1}(x_i(\xi)) + w_1(\xi) u_{2k+1}'(\xi) + \sum_{i=2}^{k+1} w_i(\xi) u_{2k+1}'(x_i(\xi))x_i'(\xi).
\end{equation}
This time we do have to include $x_{k+1}(\xi)$. The Jacobian provides the following system of equations
\begin{equation}\label{eq:exact2_derivative}
 \sum_{i=1}^{k+1} w_i' u_j(x_i) + w_1 u_j'(\xi) + \sum_{i=2}^{k+1} w_i u_j'(x_i) x_i' = 0, \qquad j = 0, \ldots, 2k.
\end{equation}


The above implementation suggestions assume that the elements of the Chebyshev set are differentiable functions. If they are not, the general theory of Chebyshev sets still applies, and so does Algorithm 2. However, we can no longer use Newton's method to solve the rootfinding problems, since it is not possible to evaluate derivatives. A suitable alternative is the simpler bisection method. Since monotonicity is guaranteed, the bisection method will converge to the single unique root.

One complication in the rootfinding process is that one can not simply evaluate $F_k(\xi)$ and $G_k(\xi)$ at arbitrary values of $\xi$. For each value of $\xi$, we do have to solve a non-linear system of equations to find the roots and weights of the canonical representation. This is the reason for our use of the word \emph{continuation}, as we proceed by small increments of $\xi$. In each step, we can start from the solution of the previous step, update to the new value of $\xi$ using the Jacobian (i.e. solving \eqref{eq:exact1_derivative} or \eqref{eq:exact2_derivative}), and use the outcome as starting values for a Newton iteration. Again, in the absence of differentiability of the basis functions, an alternative scheme would have to be devised. In our implementation, we have chosen to first approximate $F_k(\xi)$ and $G_k(\xi)$ as functions of $\xi$, using their evaluation in a sufficiently large number of samples computed via continuation. This is repeated adaptively until an accuracy threshold is reached, and the result is stored as an expansion in Chebyshev polynomials. Next, we perform the bisection algorithm on the approximating functions in order to locate the roots. The goal of this implementation is not maximal efficiency, but maximal robustness.

\section{Numerical examples}\label{s:examples}

The examples in this section have been computed using the software package \emph{GeneralizedGauss.jl}, written in the Julia programming language. The package provides a flexible framework for defining Chebyshev sets and for computing the associated quadrature rules. At the time of writing, the package is publicly available on the software repository GitHub.\footnote{\url{https://github.com/daanhb/GeneralizedGauss.jl}}

All computations in this section have been carried out in standard double precision arithmetic. However, as indicated below, for some Chebyshev sets the computation of rules with a larger number of points would require higher precision. The Julia package seamlessly allows variable-precision arithmetic with reasonable efficiency. We postpone a discussion of the accuracy of the quadrature rule to \S\ref{s:accuracy}.

\subsection{Regular Gaussian quadrature}

Polynomials up to degree $2n-1$ form a Chebyshev set, and our algorithm can be used as an alternative to the Golub-Welsch algorithm \cite{golub1969gauss,gautschi2004opq}. Note that we do not advocate this choice in practice, since Golub--Welsch is considerably more efficient for regular Gaussian quadrature (at least, once the recurrence coefficients of the associated orthogonal polynomials are available, or have been computed). Yet, the example serves well to illustrate the principles of the algorithm of this paper.

\begin{figure}
\begin{center}
 \includegraphics[width=10cm]{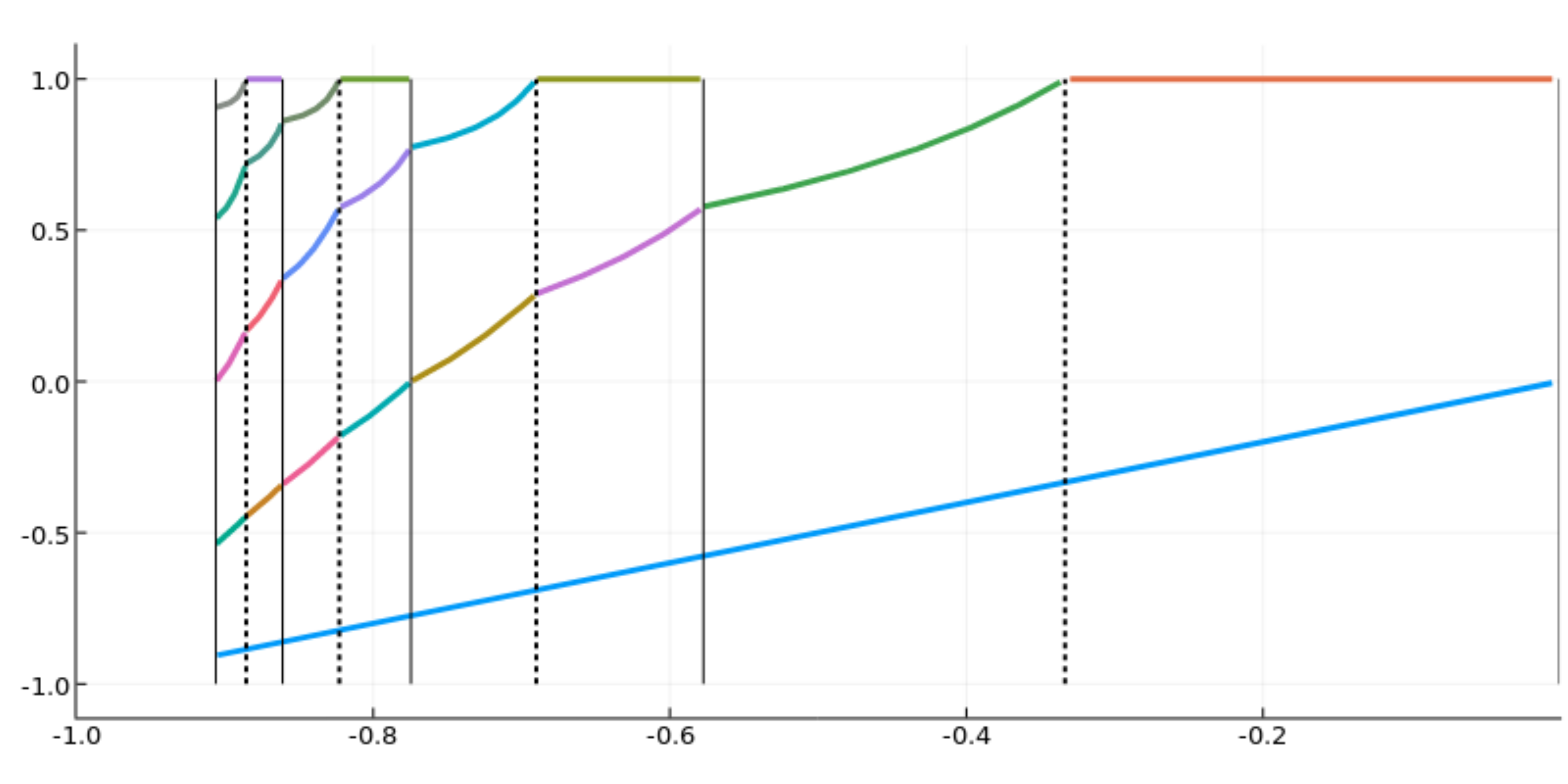}
\end{center}
\caption{Convergence for $5$-point Gauss-Legendre quadrature. The x-axis corresponds to the $\xi$ variable, which ranges from $0$ down to $-1$ in this case (from right to left). The range on the y-axis corresponds to the integration interval $[-1,1]$. The graphs display the piecewise smooth functions $x_l(\xi)$, $l=1,\ldots,5$. The vertical lines indicate the critical values of $\xi$ corresponding to lower principal representations (solid line) and upper principal representations (dotted line). Hence, the solid lines show the Gaussian quadrature rules with $1,2,3,4,5$ points respectively (from right to left), while the dotted lines show the Gauss-Lobatto rules that involve the right endpoint $1$.}\label{fig:legendre}
\end{figure}

We have computed a $5$-point Gauss-Legendre rule using the methods of this paper for the Chebyshev set
\[
 T_{2l} = \{ C_j \}_{j=0}^{2l-1},
\]
where $C_j$ are the Chebyshev polynomials of the first kind.\footnote{We avoid the standard notation $T_j$ of Chebyshev polynomials, since that notation is reserved for a Chebyshev set in this paper, following \cite{karlin1966tchebysheff}.} We remark on our choice of basis for the space of polynomials further on in \S\ref{s:accuracy}.

The behaviour of the points $x_i(\xi)$, for $i=1,\ldots,5$, is illustrated in Figure~\ref{fig:legendre}. The lower line in the figure is the first quadrature point $x_1(\xi) = \xi$. We start from a one-point rule at $\xi=0$ in the right of the figure. The point moves linearly towards $-1$. The critical values of $\xi$ are clearly visible as \emph{kinks} in the curves. At each kink, we start a new smooth continuation procedure until the next critical value of $\xi$. Every second kink, a new quadrature point enters the interval $[-1,1]$ as we look for the next lower principal representation. The vertical lines indicate the critical values of $\xi$: at each of these values, the points $x_i(\xi)$ correspond to a Gauss or Gauss-Lobatto quadrature rule on $[-1,1]$, depending on whether they constitute an upper or lower principal representation.

\subsection{Integrals with a logarithmic endpoint singularity}

An interesting application of generalized Gaussian quadrature is for integrands of the form \eqref{eq:singular_f} with a logarithmic singularity, i.e.,
\[
 f(t) = p(t) + \log(t) q(t),
\]
where $p$ and $q$ are smooth functions on $[0,1]$. If it is possible to evaluate $p$ and $q$ separately, one can write the original integral of $f$ as a sum of two integrals:
\[
 \int_0^1 w(t) f(t) {\rm d}t = \int_0^1 w(t) p(t) {\rm d}t + \int_0^1 w(t) \log(t) q(t) {\rm d}t.
\]
The first of these integrals can be evaluated using traditional Gaussian quadrature, based on polynomials that are orthogonal with respect to $w(t)$ on $[0,1]$. The second integral can be treated similarly, but using the weight function $w(t) \log(t)$ instead. The first quadrature rule requires evaluations of $p$, the second of $q$.

The value of generalized Gaussian quadrature arises when it is known that smooth functions $p$ and $q$ exist, but they can not be independently evaluated. In that case, the generalized Gaussian quadrature rule that is exact on the span of the set $\{ t^j \}_{j=0}^{l-1} \cup \{ t^j \log t \}_{j=0}^{l-1}$ leads to the simpler approximation
\[
 \int_0^1 w(t) f(t) {\rm d}t \approx \sum_{i=1}^l \lambda_i f(t_i).
\]
Here, the only requirement is a computational procedure to evaluate $f$ itself.

\begin{figure}[t]
\begin{center}
 \includegraphics[width=10cm]{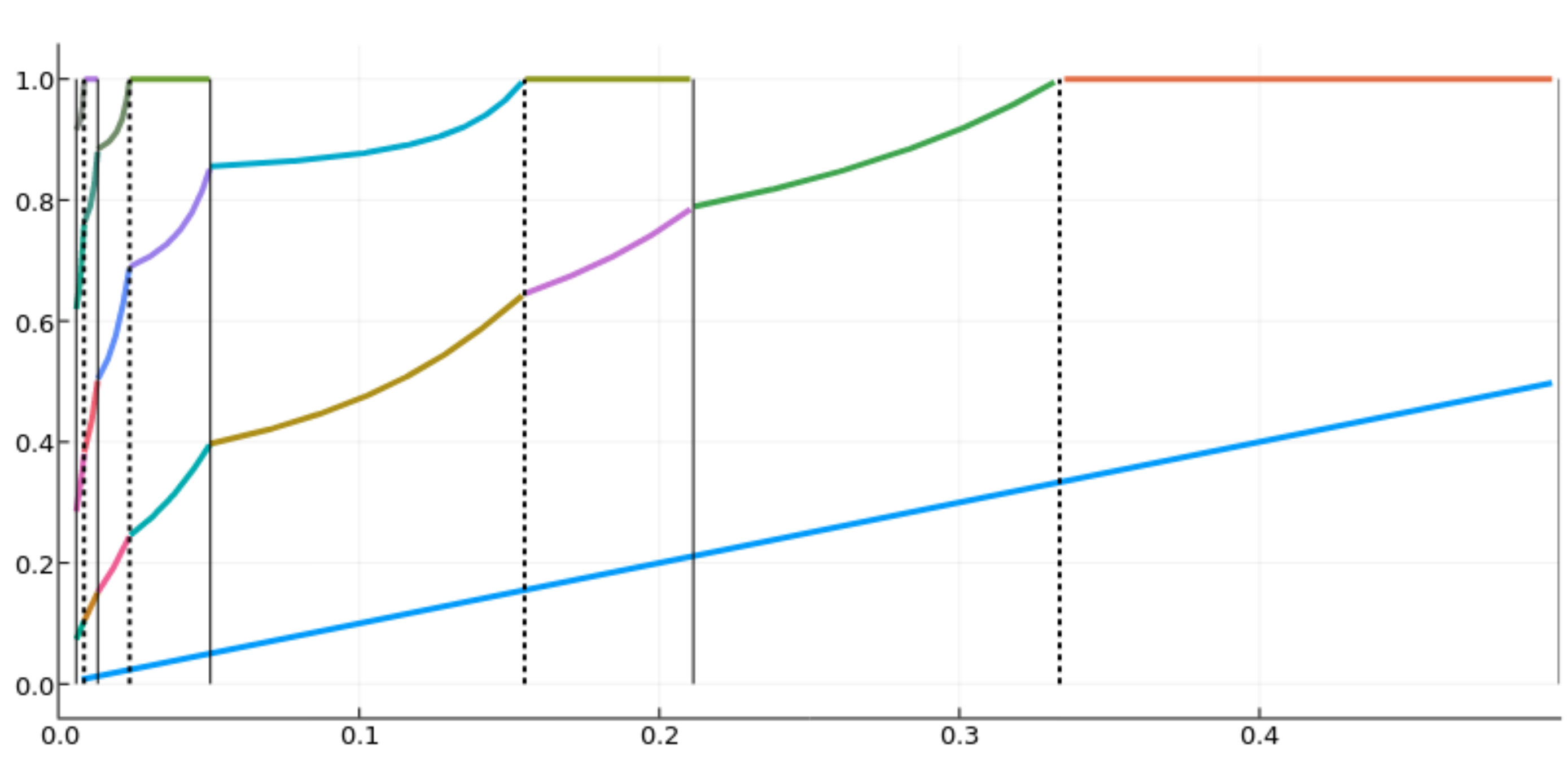}
\end{center}
\caption{Similar illustration as in Fig.~\ref{fig:legendre}, but for the set $T_{2l}$ given by \eqref{eq:log_set} with $l=5$. This set includes logarithmically singular functions, hence the $5$-point quadrature rule corresponding to the leftmost vertical line is a generalized Gaussian quadrature rule. It is well suited for integrals on $[0,1]$ with a logarithmic singularity at the left endpoint.}\label{fig:log}
\end{figure}

The convergence of this quadrature rule for logarithmically singular integrals, with corresponding error estimates, was illustrated in \cite{huybrechs2009generalizedgauss}. Here, we will only illustrate the properties of the algorithm of the current paper. We use the set
\begin{equation}\label{eq:log_set}
 T_{2l} := \{ C_j(2x-1) \}_{j=0}^{l-1} \cup \{ C_j(2x-1) \log x \}_{j=0}^{l-1}.
\end{equation}
on the interval $[0,1]$. The functions $C_j(2x-1)$ are Chebyshev polynomials as before, but scaled here to $[0,1]$.

Fig.~\ref{fig:log} is the analogue of Fig.~\ref{fig:legendre} for $l=5$. The figure shows the convergence of the algorithm (for decreasing $\xi$) to the $5$-point generalized Gaussian quadrature rule. Note that the rightmost part of figures \ref{fig:log} and \ref{fig:legendre} are in complete agreement, up to the rescaling of $[-1,1]$ to $[0,1]$. That is because the algorithm for the set $T_{2l}$ above considers first the polynomial functions $\{ C_j(2x-1) \}_{j=0}^4$, i.e. the first part of \eqref{eq:log_set}. The logarithmically singular functions are considered only afterwards. Thus, the first two quadrature rules that are computed are the regular Gauss-Legendre rules with $1$ and $2$ points respectively, scaled to $[0,1]$. A straightforward optimization would be to initiate the algorithm of this paper starting at the $\lfloor \frac{l}{2} \rfloor$-point Gauss-Legendre rule, rather than the $1$-point rule.

\begin{figure}[t]
\begin{center}
 \includegraphics[width=10cm]{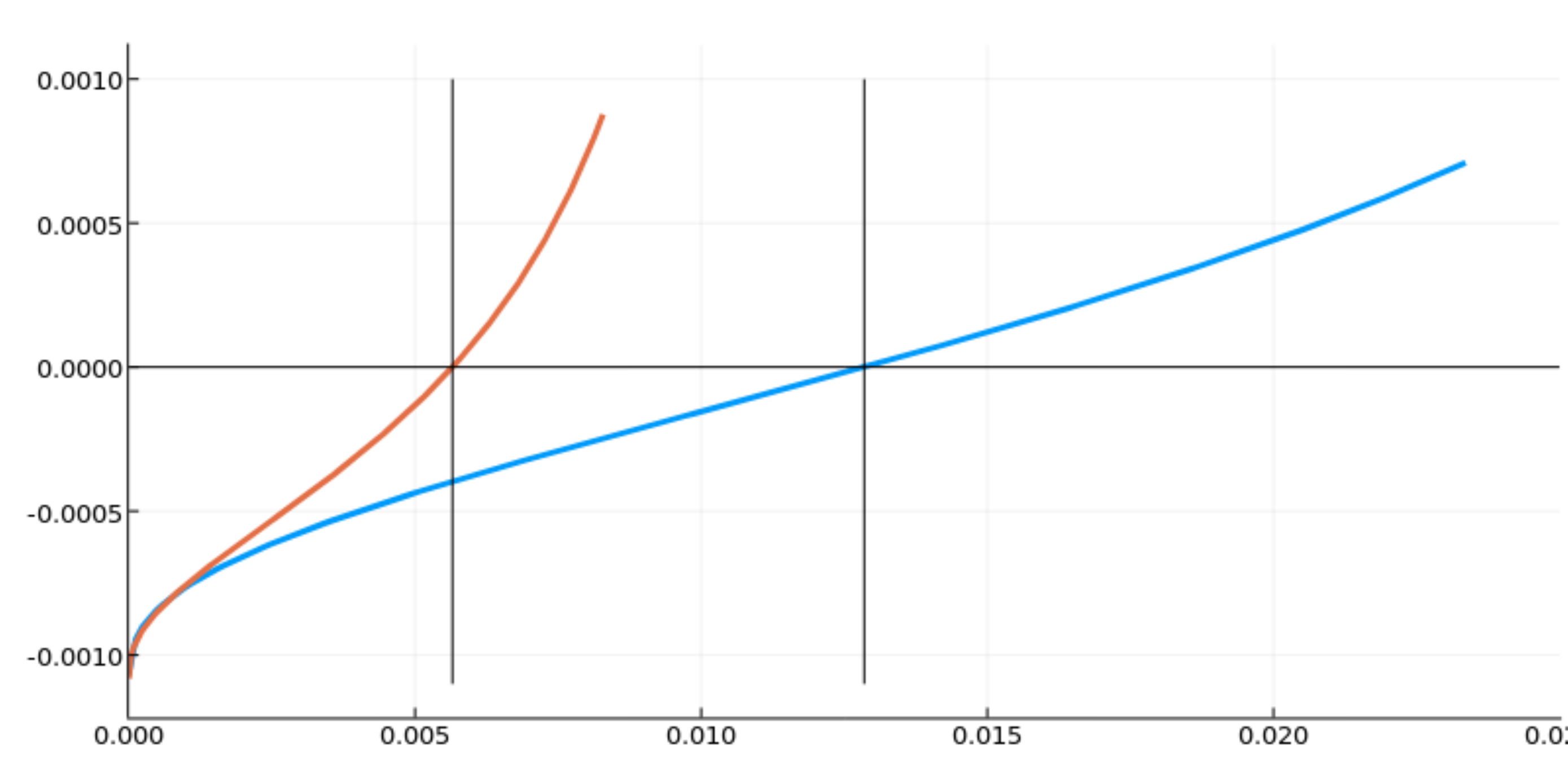}
\end{center}
\caption{Illustration of the functions $G_k$ defined in \eqref{eq:thirdproblem}. The algorithm finds the roots of these functions via continuation. The geometric theory of Chebyshev sets guarantees that $G_k$ is strictly monotonically increasing. The simplicity of this rootfinding problem for a smooth and monotonic function lies at the heart of the algorithm of this paper. Shown with vertical lines are the critical values $\xi \approx 0.00565317$ and $\xi \approx 0.012493$, which correspond to two of the lower principal representations shown in Fig.~\ref{fig:log}.}\label{fig:continuation}
\end{figure}

We conclude by illustrating the result that the functions $F_k$ and $G_k$, defined by \eqref{eq:secondproblem} and \eqref{eq:thirdproblem} respectively, are motonically increasing. Recall that the roots of these functions are the critical values of $\xi$ that are located via continuation in $\xi$. Even for the logarithmically singular set $T_{2l}$, these functions are perfectly smooth. Each has a single unique root, and the root is simple. Two such functions $G_k$ are illustrated in Fig.~\ref{fig:continuation}.

\subsection{A singular and highly oscillatory integral}

Quadrature rules for logarithmically singular functions are interesting, yet they fall short of the remarkable generality of Chebyshev sets. Several types of singular integrals arise often in boundary element methods, and several applications of generalized Gaussian quadrature rules are used in that setting in \cite{kolm2001quadrature,kolm2003quadruple}.

We consider another integral from the literature of boundary element methods that combines two difficulties: it is both singular and highly oscillatory. The model form is (see \cite{huybrechs2007highf,chandlerwilde2012actanumerica,dominguez2007hybrid} for specific examples)
\begin{equation}\label{eq:osc_int}
 I[f] =  \int_a^b f(x) H_0^{(1)}(k g_1(x)) e^{i k g_2(x)} \, {\rm d}x.
\end{equation}
Here, $H_0^{(1)}(z)$ is the Hankel function of the first kind and order zero, which has a logarithmic singularity at $z=0$ \cite{dlmf}. The functions $g_1$ and $g_2$ are phase functions, that are known explicitly in the application. The parameter $k$ corresponds to a wavenumber or frequency: if it is large, the integral above is highly oscillatory. The integral is also singular if $g_1(x) = 0$ for some $x \in [a,b]$. The value of the integral for large $k$ is mostly determined by the behaviour of $f$ near the endpoints $a$ and $b$, near any point where $g_1(x)=0$, and near so-called \emph{stationary points} where $g_1'(x)+g_2'(x)=0$ \cite{dhi2018book}.

At each of these contributing points, assuming analyticity of $f$, $g_1$ and $g_2$, the path of integration can be deformed onto the \emph{steepest descent path} along which the integrand decays exponentially and is free of oscillations \cite{huybrechs2006osc1,huybrechs2007highf,dhi2018book}. Abstracting away the details, this results in an integral of the form
\begin{equation}\label{eq:nsd_integral}
 \int_0^\infty F(t) e^{-t} \; {\rm d}t,
\end{equation}
where $F(t)$ can be evalued in terms of the function $f$ in \eqref{eq:osc_int}. This integrand potentially has both a logarithmic and a smooth part at $t=0$. Interestingly, though the Hankel function itself in \eqref{eq:osc_int} can be separated into a singular and a smooth part, neither of those parts exhibits exponential decay in the complex plane individually. It is only their combination that allows the application of the numerical method of steepest descent. Generalized Gaussian quadrature can not be avoided in this instance.

We use a Chebyshev set consisting of polynomials $\{ p_j\}_{j=0}^{l-1}$ and the logarithm times polynomials, as in the previous section, but consider integration with respect to the weight function $e^{-x}$ on the positive halfline $[0,\infty)$. For completeness, the non-linear system of equations to solve becomes
\[
\left\{
\begin{array}{ll}
\int_0^\infty p_j(x) e^{-x} {\rm d}x &= \sum_{i=1}^n \lambda_i \, p_j(x_i), \\
\int_0^\infty \log(x) p_j(x) e^{-x} {\rm d}x &= \sum_{i=1}^n \lambda_i \log(x_i) p_j(x_i),
\end{array}
\right.
\qquad j=0,1,\ldots, l-1.
\]

In order to illustrate the concept we consider a simple case where $g_1(x)=g_2(x)=x$, since in that case the steepest descent paths are fully explicit. We can write the integral \eqref{eq:osc_int} as a sum of two line integrals in the complex plane: one extends from $x=0$ upwards along the imaginary axis to infinity, the second one returns from infinity along a line parallel to the imaginary axis and ending at the other endpoint $x=1$. That is, assuming sufficient analyticity of $f$ to justify the deformation, we arrive after suitable scaling at
\begin{align*}
 I[f] &= \int_0^1 f(x) H_0^{(1)}(kx) e^{i k x} {\rm d}x \\ &= \frac{i}{2k} \int_0^\infty \left[ f\left( i \frac{t}{2k} \right)H_0^{(1)}\left( i \frac{t}{2}\right) e^{t/2} \right] e^{-t} {\rm d}t \\ &\quad - \frac{i}{2k}e^{ki} \int_0^\infty \left[ f\left( 1 + i \frac{t}{2k} \right)H_0^{(1)}\left( k + i \frac{t}{2}\right) e^{t/2} \right] e^{-t} {\rm d}t.
\end{align*}
Note that we have explicitly factored out the exponential decay of the Hankel function in the complex plane. The first integral above is weakly singular and we use the generalized Gaussian quadrature rule with $2K$ points. The second one is regular and we invoke conventional Gauss-Laguerre with $K$ points, i.e. half the number of points as for the previous integral. The results are summarized in the table below for the function $f(x) = \cos x + \sin x$. The rapid improvement with increasing $K$ should be clear, while $K$ remains fairly modest. Furthermore, the results are robust in the frequency parameter $k$, due to the steepest descent deformation.\footnote{For completeness, note that the error can be made to \emph{improve} rapidly with $k$ by using more specialized quadrature rules \cite{asheim2013transforms}. The combination of those with the generalized Gaussian quadrature rules in this paper is a promising topic of further research.}

\begin{table}[ht]
\caption{Absolute error of the approximation of \eqref{eq:osc_int} with $f(x) = \cos(x)+\sin(x)$, $g_1(x)=g_2(x)=x$, using steepest descent deformation followed by generalized Gaussian quadrature with $3K$ points in total.}\label{t:error}
\footnotesize
\begin{center}
\begin{tabular}{|c|cccc|}\hline
 $K \setminus k$ & $10$ & $20$ & $30$ & $40$ \\ \hline
 $1$ & $2.2e-4$ & $1.1e-4$ & $7.5e-5$ & $5.6e-5$ \\
 $2$ & $1.2e-6$ & $5.7e-7$ & $3.8e-7$ & $2.8e-7$ \\
 $3$ & $6.0e-9$ & $2.9e-9$ & $1.9e-9$ & $1.5e-9$ \\
 $4$ & $2.1e-11$ & $9.2e-12$ & $6.0e-12$ & $4.4e-12$ \\\hline
\end{tabular}
\end{center}
\end{table}


Though integral \eqref{eq:osc_int} is widely perceived as being very challening, our experiment shows that its complications are limited to the preparatory work. This work involves locating the critical points of the integral, in this case points where $g_1(x)=0$ and points where $g_1'(x)+g_2'(x)=0$, and finding the corresponding steepest descent paths. This can be done either analytically, as in the simple example above, or numerically. Afterwards, the evaluation of the resulting integrals of the form \eqref{eq:nsd_integral} is nearly trivial, owing to the flexibility of Chebyshev sets compared to polynomials.

%
%
%
%

\section{Accuracy of the quadrature rules}
\label{s:accuracy}

What is the numerical accuracy of the computed quadrature rules? Not of their application in the approximation of integrals -- what is the accuracy of the points and weights themselves? This is a relevant question which we do not fully settle in the current paper, as it leads to a rich topic of its own. Indeed, in the same way as Chebyshev sets potentially give rise to unconventional quadrature rules, they also give rise to unconventional approximation schemes. We do touch upon the main issues at hand.

\subsection{A source of errors: the interpolation problem}

The algorithm essentially performs rootfinding on the functions $F_k$ and $G_k$, defined by \eqref{eq:secondproblem} and \eqref{eq:thirdproblem}. In our examples we have chosen the simple bisection method combined with continuation in $\xi$, and as part of the latter we have to solve the linear systems \eqref{eq:exact1_derivative} and \eqref{eq:exact2_derivative}. These linear systems are the Jacobians of the non-linear systems of equations \eqref{eq:exact1} and \eqref{eq:exact2}, which describe the exactness of the quadrature rules.

The accuracy with which we can solve these linear systems depends on their condition number. In turn, this depends on the stability properties of the functions $u_j$ in the set $T_{2l} = \{u_j\}_{j=0}^{2l-1}$ as a basis for their span. Indeed, while the quadrature rule itself only depends on the \emph{function space} ${\mathcal T}_{2l} = \SPAN T_{2l}$, the condition numbers of \eqref{eq:exact1_derivative}--\eqref{eq:exact2_derivative} clearly depend on the basis that is chosen for that space. While a Chebyshev set $T_{n}$ on $[a,b]$ has the mathematical property that the interpolation problem is uniquely solvable for any set of $n$ distinct points in $[a,b]$, there is no guarantee that all these interpolation problems are well-conditioned. Indeed, we could have used the monomial basis in the previous section, rather than Chebyshev polynomials, and that would certainly have led to poor numerical results.

\emph{Note that this observation about the stability of the basis $\{ u_j \}$ is true for any scheme that solves the non-linear systems numerically using their Jacobians, not just the scheme of this paper.}

In practice, though there is no strict equivalence, fortunately there is a close agreement between points that are suitable for interpolation, and points of a Gaussian quadrature rule. For example, it is well-known that polynomial interpolation is best performed in roots of orthogonal polynomials, i.e., in the Gaussian quadrature points \cite{trefethen2012atap}. Furthermore, by construction, all the point sets that arise in our algorithm correspond to a quadrature rule with high order of exactness. This is true even while performing the continuation steps, as all intermediate rules correspond at least to canonical representations.

What is a good basis? One possibility beyond orthogonal bases is to consider a Riesz basis \cite{christensen2003frames}. The set $\{ u_j \}_{j=0}^n$ is a Riesz basis for its span ${\mathcal T}_n$ if there exist constants $A_n,B_n > 0$ such that for every finite scalar sequence $\{c_k\}$ one has
\begin{equation}\label{eq:Riesz}
 A_n \sum_{k=0}^{n} c_k^2  \leq \left\Vert \sum_{j=0}^n c_j  u_j \right\Vert^2 \leq  B_n \sum_{k=0}^{n} c_k^2.
\end{equation}
Stability of a Riesz basis (in the terminology of approximation theory) hinges on the ratio $B_n/A_n$, which is always greater than or equal to $1$. The closer this ratio is to $1$, the better. Of interest in the current setting is the possible growth of $B_n / A_n$ for increasing $n$, which is to be avoided. Ideally, the ratio is uniformly bounded in $n$.

Yet, this avenue requires at least a norm on ${\mathcal T}_n$, which we have thus far managed to avoid. Moreover, many interesting Chebyshev sets are \emph{not} stable Riesz bases for increasing $n$.

\subsection{Chebyshev sets and truncated frames}

The second example in the previous section, the set of polynomials and $\log$ times polynomials given by \eqref{eq:log_set}, gives rise to an exponentially growing ratio $B_n / A_n$. Rather than as an ill-conditioned basis, it is best seen as the truncation of an infinite \emph{frame} \cite{christensen2003frames,huybrechs2016tw674}. A set is a frame for a separable Hilbert space $H$ if it is dense in $H$, though possibly redundant, while satisfying the so-called frame condition
\begin{equation}\label{eq:framecondition}
 A \Vert f \Vert \leq \sum_{j=1}^\infty \langle f, u_j \rangle^2 \leq B \Vert f \Vert, \qquad \forall f \in H.
\end{equation}

Related to our second example, the infinite set
\begin{equation}\label{eq:frame}
 \Phi := \{ T_j(2x-1) \}_{j=0}^\infty \cup \{ T_j(2x-1) \log x \}_{j=0}^\infty.
\end{equation}
is a frame for $L^2([\delta,1])$ for any $\delta > 0$. It is clearly complete, since it includes all polynomials which are dense in $L^2$ on any interval. However, it is also redundant, since the polynomials alone are already sufficient for denseness. Yet, though polynomials are dense in $L^2$, it is clear that one requires a large degree in order to approximate a (nearly) singular function to any satisfiable accuracy. In contrast, accurate approximations in the span of a truncation of the frame \eqref{eq:frame} clearly require only a modest number of degrees of freedom.

Unfortunately, the frame condition \eqref{eq:framecondition} for the infinite frame does not give rise to conditions of the form \eqref{eq:Riesz} with bounded ratios $B_n / A_n$ after truncation: truncated frames tend to lead to extremely ill-conditioned linear systems in approximation problems. This effect is described and explained in detail in \cite{huybrechs2016tw674}. This remains true even if the truncated frame (of which \eqref{eq:log_set} is an example) contains only linearly independent functions and if it is a finite Chebyshev set.

It is shown in \cite{huybrechs2016tw674} that, in spite of potentiall ill-conditioning, approximations in truncated frames can be computed accurately in standard finite precision using conventional regularization techniques, in a least squares sense. This does not immediately apply to the current algorithm: for the time being, example 2 in this paper requires high precision arithmetic for large $n$. It does appear that the theory of frames may suggest an alternative basis for the same function space. Furthermore, initial experiments indicate that generalized Gaussian quadrature points yield suitable points for \emph{interpolation} in a truncated frame, as opposed to the least squares approximation that is suggested in \cite{huybrechs2016tw674}. This is a promising avenue of ongoing and future research.

\section{Concluding remarks}

To the best of our knowledge, the algorithm presented in this paper is the first one to compute generalized Gaussian quadrature rules for complete Chebyshev sets with guaranteed success. It enables novel quadrature approximation for many kinds of integrals that are not well evaluated with polynomial-based quadrature. Yet, it is important to appreciate the limitations. First, the algorithm is not as general as the theory of Chebyshev sets actually allows. Indeed, several Chebyshev sets of practical interest are not complete Chebyshev sets. Some other interesting sets, including highly oscillatory functions for example, are not Chebyshev sets at all. Second, the stability properties of numerical approximations using general Chebyshev sets are not as well understood as they are for polynomial approximations. It may often be possible to identify useful Chebyshev sets, the span of which would contain good approximations to a given integrand under investigation. However, it may require additional work to compute the corresponding quadrature rule in a numerically stable way in finite precision arithmetic. Still, the quest for a good basis in a particular function space is decoupled from an algorithm to compute the corresponding quadrature rule, and in this paper we have largely addressed at least the latter part.

\bibliographystyle{abbrv}
\bibliography{references}

\end{document}